# ON THE NUMBER OF SUPPORT POINTS OF MAXIMIN AND BAYESIAN OPTIMAL DESIGNS[1]

BY DIETRICH BRAESS AND HOLGER DETTE

*Ruhr-Universität Bochum*

We consider maximin and Bayesian $D$-optimal designs for nonlinear regression models. The maximin criterion requires the specification of a region for the nonlinear parameters in the model, while the Bayesian optimality criterion assumes that a prior for these parameters is available. On interval parameter spaces, it was observed empirically by many authors that an increase of uncertainty in the prior information (i.e., a larger range for the parameter space in the maximin criterion or a larger variance of the prior in the Bayesian criterion) yields a larger number of support points of the corresponding optimal designs. In this paper, we present analytic tools which are used to prove this phenomenon in concrete situations. The proposed methodology can be used to explain many empirically observed results in the literature. Moreover, it explains why maximin $D$-optimal designs are usually supported at more points than Bayesian $D$-optimal designs.

**1. Optimal designs for nonlinear models.** Consider the common problem of nonlinear experimental design where the scalar response variable, say $Y$, is distributed as a member of the exponential family with

$$(1.1) \qquad E[Y|x] = \eta(x, \theta),$$

$\theta \in \mathbb{R}^m$ is the unknown parameter, $x$ denotes the explanatory variable that varies in a compact space, say $\mathcal{X}$, and $\eta$ is a given function. We assume that observations under different experimental conditions are independent and denote the Fisher information matrix for the parameter $\theta$ at the point $x$ by

$$(1.2) \qquad I(x, \theta) = \frac{1}{\text{Var}[Y|x]} \left( \frac{\partial \eta}{\partial \theta}(x, \theta) \right) \left( \frac{\partial \eta}{\partial \theta}(x, \theta) \right)^T \in \mathbb{R}^{m \times m};$$

Received September 2005; revised June 2006.
[1]Supported by the Deutsche Forschungsgemeinschaft (SFB 475, "Komplexitätsreduktion in multivariaten Datenstrukturen").

see [7]. Throughout this paper, we assume the continuous differentiability of the function $\eta$ with respect to $\theta$ and the existence of the conditional variance.

An approximate design $\xi$ for this model is a probability measure on the design space $\mathcal{X}$ with finite support $x_1, \ldots, x_n$ and weights $w_1, \ldots, w_n$ representing the relative proportions of total observations taken at the corresponding design points; see, for example, [11]. The information matrix of a design $\xi$ is defined by

$$(1.3) \qquad M(\xi, \theta) = \int_{\mathcal{X}} I(x, \theta) \, d\xi(x),$$

and a local optimal design maximizes a given function of the matrix $M(\xi, \theta)$; see [18]. We consider the $D$-optimality criterion $\log |M(\xi, \theta)|$ for discriminating among competing designs. In general, such a design depends on the unknown parameter $\theta$, which must be specified for its implementation. Local optimality criteria have been criticized by numerous authors because the resulting optimal designs can be highly inefficient within the true model setting if the unknown parameters are misspecified.

A more robust approach has been achieved in practice by using the concepts of Bayesian and maximin optimality, since additional information on the uncertainty in those parameters can be incorporated. A priori knowledge of the experimenter can be modeled mathematically as follows. Assume that $\theta \in \Theta$, where $\Theta \subset \mathbb{R}^m$ denotes a set and let $\pi$ denote a probability measure on $\Theta$. A design is called *Bayesian $D$-optimal* (with respect to a given prior $\pi$ on $\Theta$) if it maximizes the function

$$(1.4) \qquad \int_{\Theta} \log |M(\xi, \theta)| \pi(d\theta).$$

Bayesian $D$-optimal designs have been studied by Chaloner and Larntz [2], Pronzato and Walter [17], Mukhopadhyay and Haines [15], Dette and Neugebauer [5, 6] and many others.

In some circumstances, it is difficult for the experimenter to specify a prior on the parameter space $\Theta$. Therefore, several authors have proposed *standardized maximin $D$-optimal designs*, that is, designs which maximize

$$(1.5) \qquad \min\left\{ \frac{|M(\xi, \theta)|}{|M(\xi[\theta], \theta)|} \Big| \theta \in \Theta \right\},$$

where $\xi[\theta]$ denotes the local $D$-optimal design for fixed $\theta$; see, for example, [3, 9, 16] or [4]. The criterion (1.5) does not compare the quantities $|M(\xi, \theta)|$ directly, but rather with respect to the values that could be obtained if $\theta$, and, as a consequence, the local $D$-optimal design, were known. Bayesian and standardized maximin $D$-optimal designs can only be given explicitly in rare circumstances. Moreover, optimal designs often require more than $m$ support points unless the parameter space $\Theta$ is "sufficiently narrow" or the prior in the Bayesian criterion puts "most of its mass at a small" subset of $\Theta$. For interval parameter spaces, it was observed empirically that the number of support points increases if less knowledge about $\theta$ is incorporated in the optimality criteria (see [2] and [4], among others).

In the present paper, we will provide analytic tools for making a rigorous decision as to whether the number of support points is unbounded if the a priori information on an interval parameter space is repeatedly diminished. It turns out that the answer to this question depends in a complicated way on the given nonlinear model. Therefore, we will mainly concentrate on the case where only one parameter, say $\beta$, enters nonlinearly in the optimality criterion and will mention possible extensions to the general case briefly in the Appendix. For the maximin criterion, we quantify uncertainty by considering an interval with increasing length, while for the Bayesian case, we use the properties of the prior—this includes a uniform prior on an interval with increasing length. We establish sufficient conditions on the nonlinear model such that increasing uncertainty about the nonlinear parameter leads to an arbitrarily large number of support points of Bayesian and standardized maximin $D$-optimal designs. In fact, the tools in this paper are applicable to all models known to us.

The conditions are more restrictive for the Bayesian $D$-optimal design than for the standardized maximin $D$-optimal one. This explains why standardized maximin $D$-optimal designs are usually supported at more points than Bayesian $D$-optimal designs. In particular, in the case of Bayesian $D$-optimality, the number of support points may increase so slowly that it is almost impossible to decide by numerical computations whether the number is asymptotically bounded or not.

In Section 3 we discuss standardized maximin optimal designs, while Section 4 deals with the Bayesian case. We illustrate our approach for models with one, two and three parameters. The proposed methodology is a general one, but the technical difficulties for the verification of the conditions differ in each scenario and increase with the dimension. Finally, some conclusions are given in Section 5, while all technical details are deferred to the Appendix. For further examples, see the technical report of Braess and Dette [1].

**2. Preliminaries.** For the sake of simplicity, we assume that the local $D$-optimal design depends only on one component of the parameter $\theta \in \mathbb{R}^m$. The general situation can be obtained by a straightforward generalization which is briefly indicated in Appendix A.3. We denote this component by $\beta$ and the corresponding design by $\xi[\beta]$. Consequently, we will reflect only this dependence in our notation and the optimality criteria in (1.4) and (1.5) are represented by the functions

$$\Psi(\xi) = \int_{\mathcal{B}} \log |M(\xi,\beta)| \pi(d\beta), \tag{2.1}$$

$$\Phi(\xi) = \min\left\{ \frac{|M(\xi,\beta)|}{|M(\xi[\beta],\beta)|} \,\Big|\, \beta \in \mathcal{B} \right\}, \tag{2.2}$$

respectively. Here, $M(\xi,\beta)$ is the information matrix (1.3) in the nonlinear model, $\mathcal{B} = [\beta_{\min}, \beta_{\max}]$ represents the prior knowledge about the location of the unknown parameter $\beta$ and $\pi$ denotes a prior on $\mathcal{B}$.

Let $\xi$ be a design on $\mathcal{X}$ with masses $w_k$ at the support points $x_k$ ($k = 1, \ldots, n$). The information matrix of $\xi$ is then given by

$$M(\xi,\beta) = \sum_{k=1}^{n} w_k I(x_k,\beta).$$

Throughout this paper,

$$Q(\beta,\tilde{\beta}) = \frac{|M(\xi[\tilde{\beta}],\beta)|}{|M(\xi[\beta],\beta)|} \tag{2.3}$$

quantifies the loss of information if $\beta$ is the "true" unknown parameter, while the experimenter uses the local $D$-optimal design for the (wrong) guess $\tilde{\beta}$. Note that $(Q(\beta,\tilde{\beta}))^{1/m}$ is the $D$-efficiency of the design $\xi[\tilde{\beta}]$. We will derive sufficient conditions such that the number of support points of the optimal designs with respect to the criteria (2.1) and (2.2) exceeds any given number if the amount of prior information is decreased.

Our first definition quantifies the loss of efficiency caused by an application of a local $D$-optimal design based on a misspecified parameter.

DEFINITION 2.1.  Let $\ell : \mathcal{B} \to \mathbb{R}$ be a nondecreasing continuous function. The function $Q$ defined in (2.3) is said to be *uniformly decreasing with respect to the scale $\ell$* if the following two conditions hold:

(i) for all $\beta, \tilde{\beta} \in \mathcal{B}$, the inequality

$$Q(\beta,\tilde{\beta}) \leq \varphi(\ell(\beta) - \ell(\tilde{\beta})) \tag{2.4}$$

holds, where $\varphi$ is a real-valued function whose decay [i.e., $\varphi(z) \to 0$] for $z \to \infty$ will be sufficiently fast, as specified later for each case under consideration;

(ii) there is a positive constant $\lambda > 0$ such that

$$Q(\beta,\tilde{\beta}) \geq \tfrac{1}{2} \quad \text{whenever } |\ell(\beta) - \ell(\tilde{\beta})| \leq \lambda. \tag{2.5}$$

There is some heuristic explanation for the two conditions in Definition 2.1. The quantity

$$d(\beta,\tilde{\beta}) = |\ell(\beta) - \ell(\tilde{\beta})| \tag{2.6}$$

can be considered a distance between the parameters $\beta$ and $\tilde{\beta}$. This interpretation is also useful for the extension to models with more than one nonlinear parameter. On one hand, condition (2.4) requires that the efficiency decrease sufficiently fast if the parameter is misspecified. On the other hand, condition (2.5) guarantees that the efficiency cannot become small if the parameter is only slightly misspecified. Since these conditions are very natural, they are satisfied by most of the commonly used nonlinear models with $\ell(\beta) = \beta$ or $\ell(\beta) = \log \beta$. In fact, we are not aware of any model where conditions (2.4) and (2.5) are not satisfied.

A final assumption is required for our main results. Roughly speaking, it guarantees that points in the design space $\mathcal{X}$ which are not in the support of any local $D$-optimal design $\xi[\beta]$ can be disregarded for the construction of the optimal designs. To be precise, we represent the Fisher information as

$$(2.7) \qquad I(x_1, \ldots, x_m, \beta) = f(x, \beta) f^T(x, \beta),$$

where $f(x, \beta) = (f_1(x, \beta), \ldots, f_m(x, \beta))^T \in \mathbb{R}^m$, and introduce, for $x_1, \ldots, x_m \in \mathcal{X}$, the determinant

$$(2.8) \qquad I_m(x_1, \ldots, x_m, \beta) = \begin{vmatrix} f_1(x_1, \beta) & \cdots & f_1(x_m, \beta) \\ \vdots & \vdots & \vdots \\ f_m(x_1, \beta) & \cdots & f_m(x_m, \beta) \end{vmatrix}^2.$$

We assume that there exists a constant, say $c_0$, such that for any $x = (x_1, \ldots, x_m)^T \in \mathcal{X}^m$, there exist local $D$-optimal designs $\xi[\tilde{\beta}^{(1)}], \ldots, \xi[\tilde{\beta}^{(m)}]$ with $\tilde{\beta}^{(1)}, \ldots, \tilde{\beta}^{(m)} \in \mathcal{B}$, such that

$$(2.9) \qquad |I_m(x_1, \ldots, x_m, \beta)| \leq c_0 \sum_{j=1}^m |M(\xi[\tilde{\beta}^{(j)}], \beta)| \qquad \forall \beta \in \mathcal{B}.$$

Additionally, we define $m_\eta \geq 0$ as the number of points which appear as support points of any local $D$-optimal design in the nonlinear model and we assume $m > m_\eta$. For example, if $\eta(x, \theta) = \theta_0 + \theta_1 e^{-\theta_2 x}$ and $x \in [0, 1]$, it follows from [8] that any local $D$-optimal design has three support points, including the points 0 and 1. For this model, we have $m = 3$ and $m_\eta = 2$.

Finally, we note that assumption (2.9) is obviously satisfied in models where the local $D$-optimal designs are supported at a minimal number of $m$ points. In examples where the local $D$-optimal designs are supported at more than $m$ points, the condition has to be checked in each situation. However, thus far, the authors have not found a case where (2.9) is not satisfied.

**3. Standardized maximin $D$-optimal designs.** The following result shows for nonlinear models that the number of support points of the standardized maximin $D$-optimal design can become arbitrarily large under the assumptions stated in Section 2.

THEOREM 3.1. *Assume that $Q$ is uniformly decreasing with respect to the scale $\ell$ in the sense of Definition* 2.1, *where*

$$\varphi(z) \leq c_1 |z|^{-\gamma} \qquad \text{with } c_1 > 0, \ \gamma > m, \tag{3.1}$$

*and that* (2.9) *is satisfied. Let $N \in \mathbb{N}$ be given. If $\ell(\beta_{\max}) - \ell(\beta_{\min})$ is sufficiently large, then the standardized maximin D-optimal design with respect to the interval $\mathcal{B} = [\beta_{\min}, \beta_{\max}]$ is supported at more than $N$ points.*

*Moreover, if the local D-optimal designs are minimally supported and the support contains only $m$ points, then the condition $\gamma > m$ in* (3.1) *can be replaced by $\gamma > m - m_\eta$.*

The technical proof is deferred to the Appendix. The main idea is to use condition (3.1) to show that the value of the standardized maximin $D$-optimality criterion of any $N$-point design can be estimated by $O(B^{-\gamma})$, where $B = \ell(\beta_{\max}) - \ell(\beta_{\min})$. In a second-step condition, (2.5) is used to construct a design with $n > N$ support points, for which the value of the criterion is at least of the order $O(B^{-m})$. Since $\gamma > m$, it follows, for sufficiently large $B$, that no design with $N$ support points can be standardized maximin $D$-optimal.

In the following we will illustrate the application of Theorem 3.1 in the cases $m = 1, 2, 3$. The technical difficulties for the verification of the sufficient conditions increase with the dimension of the Fisher information.

EXAMPLE 3.2. Consider the one-dimensional exponential growth model with normally distributed homoscedastic errors and

$$\eta(x, \beta) = e^{-\beta x}, \qquad \beta \in [1, B], x \in [0, 1]. \tag{3.2}$$

The Fisher information of the parameter $\beta$ (up to a constant which does not affect the optimal design problem) is

$$I(x, \beta) = x^2 e^{-2\beta x}. \tag{3.3}$$

The local $D$-optimal design is a one-point design supported at the point $x[\beta] = 1/\beta$ and it follows from

$$M(\xi[\beta], \beta) = I(x[\beta], \beta) = (e\beta)^{-2} \tag{3.4}$$

that the function $Q$ in (2.3) is given by

$$Q(\beta, \tilde{\beta}) = \left(\frac{\beta}{\tilde{\beta}} e^{1-\beta/\tilde{\beta}}\right)^2. \tag{3.5}$$

Hence, $Q(\beta, \tilde{\beta}) = \psi(\frac{\beta}{\tilde{\beta}})$, where

$$\psi(y) = (ye^{1-y})^2 \leq \begin{cases} e^2 y^2, & \text{if } y \leq 1, \\ 3y^{-2}, & \text{if } y > 1. \end{cases}$$

TABLE 1
*Standardized maximin D-optimal designs for the exponential regression model (3.2) on the interval $[0,1]$ with respect to various parameter spaces $[1, B]$. First row: support points; second row: weights*

| B | 10 | 40 | 50 | 100 | 200 |
|---|---|---|---|---|---|
| $x_1$ | 0.142 | 0.037 | 0.028 | 0.014 | 0.007 |
| $w_1$ | 0.553 | 0.414 | 0.379 | 0.336 | 0.306 |
| $x_2$ | 0.771 | 0.193 | 0.131 | 0.064 | 0.034 |
| $w_2$ | 0.447 | 0.272 | 0.221 | 0.193 | 0.182 |
| $x_3$ | | 0.772 | 0.374 | 0.156 | 0.101 |
| $w_3$ | | 0.314 | 0.170 | 0.093 | 0.147 |
| $x_4$ | | | 0.972 | 0.287 | 0.250 |
| $w_4$ | | | 0.230 | 0.137 | 0.089 |
| $x_5$ | | | | 0.838 | 0.326 |
| $w_5$ | | | | 0.241 | 0.066 |
| $x_6$ | | | | | 0.856 |
| $w_6$ | | | | | 0.210 |

We choose $\ell(\beta) = \log \beta$ and (2.4) holds with

$$\varphi(z) \leq e^2 e^{-2|z|}. \tag{3.6}$$

The decay is even faster than required in (3.1). Moreover, $\psi(z) \geq \frac{1}{2}$ if $\frac{1}{2} \leq z \leq 2$, which proves property (ii) in Definition 2.1. Finally, we verify property (2.9) for $m = 1$. Consider a point, say $x_0$. If $1 \geq x_0 \geq 1/B$, then we have $\delta_{x_0} = \xi[1/x_0]$ for the Dirac measure at the point $x_0$ and there is, in fact, equality in (2.9) with $\tilde\beta = 1/x_0$. On the other hand, if $x_0 < 1/B$, then the Dirac measure $\delta_{x_0}$ is not local $D$-optimal for any $\beta \in [1, B]$ and $x_0\beta < \frac{\beta}{B} \leq 1$ for all $\beta \in [1, B]$. Since the function $z \mapsto z^2 e^{-2z}$ is increasing on the interval $[0, 1]$, it follows that

$$\beta^2 I_1(x_0, \beta) = \beta^2 M(\delta_{x_0}, \beta)$$
$$\leq \beta^2 M(\delta_{1/B}, \beta) = \beta^2 M(\xi[B], \beta) \quad \text{for all } \beta \in [1, B],$$

which shows (2.9). Similarly, if $x_0 > 1$, we obtain $\tilde\beta = 1$. Therefore, condition (2.9) is satisfied with

$$\tilde\beta = \max\{1, \min\{B, 1/x_0\}\}.$$

By Theorem 3.1, the number of support points of the standardized maximin $D$-optimal design for the regression model (3.2) becomes arbitrarily large with increasing parameter $B \to \infty$.

Numerical results in Table 1 illustrate this fact. We have calculated the standardized maximin $D$-optimal designs using Matlab for various parameter spaces $\mathcal{B} = [1, B]$. The optimality of the calculated designs was checked by the equivalence theorem of Wong [19].

EXAMPLE 3.3. Consider the exponential growth model with normally distributed homoscedastic errors and

$$(3.7) \qquad \eta(x, \alpha, \beta) = \alpha + e^{-\beta x}, \qquad x \in [0, 1], \ \beta \in [1, \beta_{\max}],$$

which is used for analyzing the growth of crops; see [13]. The Fisher information matrix (up to a constant which does not change the optimal design problem) for the parameter $\theta = (\alpha, \beta)$ is given by

$$I(x, \beta) = \begin{pmatrix} 1 & -xe^{-\beta x} \\ -xe^{-\beta x} & x^2 e^{-2\beta x} \end{pmatrix}$$

and the function $I_2$ defined in (2.8) is determined as

$$I_2(x_1, x_2, \beta) = (x_1 e^{-\beta x_1} - x_2 e^{-\beta x_2})^2.$$

If $x_1 < x_2$, we can always increase $I_2(x_1, x_2, \beta)$ by using $x_1 = 0$. It follows from [6] that the local $D$-optimal design puts equal masses at the points $x_1 = 0$ and $x_2 = \frac{1}{\beta}$ with corresponding determinant

$$|M(\xi[\beta], \beta)| = \frac{1}{4(e\beta)^2};$$

see also [8]. Therefore, we obtain $m = 2$, $m_\eta = 1$, $m - m_\eta = 1$ for the quantities in the second part of Theorem 3.1. From $I_2(0, 1/\tilde{\beta}, \beta) = \tilde{\beta}^{-2} e^{-2\beta/\tilde{\beta}}$, it follows that $Q(\beta, \tilde{\beta}) = \psi(\beta/\tilde{\beta})$ with the same function $\psi(z) = z^2 e^{2(1-z)}$ as in Example 3.2. This shows that (2.4) and (3.1) are satisfied. Obviously,

$$I_2(x_1, x_2, \beta) \le (x_1 e^{-\beta x_1})^2 + (x_2 e^{-\beta x_2})^2 = I_2(0, x_1, \beta) + I_2(0, x_2, \beta)$$

and we conclude, as in the previous example, that the remaining assumption (2.9) of Theorem 3.1 is also satisfied. Therefore, the number of support points of the standardized maximin $D$-optimal design in the exponential growth model (3.7) becomes arbitrarily large if $\beta_{\max} \to \infty$.

EXAMPLE 3.4. To illustrate how the technical difficulties increase with the dimension, we consider the exponential growth model with normally distributed homoscedastic errors and

$$(3.8) \qquad \eta(x, \alpha_1, \alpha_2, \beta) = \alpha_1 + \alpha_2 e^{-\beta x}, \qquad x \in [0, 1], \ \beta \in [1, \beta_{\max}].$$

In particular, we obtain for the determinant in (2.8)

$$\begin{aligned}
I_3(x_1, x_2, x_3, \beta) &= [x_1 e^{-\beta x_1}(e^{-\beta x_3} - e^{-\beta x_2}) \\
&\quad + x_2 e^{-\beta x_2}(e^{-\beta x_1} - e^{-\beta x_3}) + x_3 e^{-\beta x_3}(e^{-\beta x_2} - e^{-\beta x_1})]^2 \\
&= H^2(x_1, x_2, x_3),
\end{aligned}$$

where the last line defines the function $H(x_1, x_2, x_3)$. Han and Chaloner [8] showed that local $D$-optimal designs for the exponential regression model (3.8) have three support points and that $\bar{x} = 0$ and $\hat{x} = 1$ are in fact common support points of all local $D$-optimal designs, that is, $m = 3$, $m_\eta = 2$, $m - m_\eta = 1$. An alternative derivation of the latter fact will be given below. For the points $\bar{x} = 0$ and $\hat{x} = 1$, the determinant reduces to

$$(3.9) \quad I_3(0, x, 1, \beta) = [xe^{-\beta x}(1 - e^{-\beta}) - e^{-\beta}(1 - e^{-\beta x})]^2 = H^2(0, x, 1, \beta).$$

Note that sums of cyclic products such as $a(b - c) + b(c - a) + c(a - b)$ vanish. Therefore, we obtain the useful representation

$$\begin{aligned}
H(x_1, x_2, x_3, \beta) \\
&= \frac{e^{-\beta x_3} - e^{-\beta x_2}}{1 - e^{-\beta}}[x_1 e^{-\beta x_1}(1 - e^{-\beta}) - e^{-\beta}(1 - e^{-\beta x_1})] \\
&\quad + \frac{e^{-\beta x_1} - e^{-\beta x_3}}{1 - e^{-\beta}}[x_2 e^{-\beta x_2}(1 - e^{-\beta}) - e^{-\beta}(1 - e^{-\beta x_2})] \\
&\quad + \frac{e^{-\beta x_2} - e^{-\beta x_1}}{1 - e^{-\beta}}[x_3 e^{-\beta x_3}(1 - e^{-\beta}) - e^{-\beta}(1 - e^{-\beta x_3})] \\
&= \sum_{k=1}^{3} a_k H(0, x_k, 1, \beta),
\end{aligned} \quad (3.10)$$

with coefficients $a_k$ satisfying $|a_k| \leq 1$. If the support points are ordered, that is,

$$(3.11) \qquad 0 \leq x_1 < x_2 < x_3 \leq 1,$$

then $a_1, a_3 < 0$ and $a_2 > 0$. The estimation of the quantity $I_3$ heavily depends on the knowledge of the signs of the function $H$. Obviously, $H(0, x, 1, \beta)$ vanishes at $x = 0$ and $x = 1$. Moreover, the derivative at $x = 0$ is positive. Since exponential sums of the form $a_1 + (a_2 + a_3 x)e^{-\beta x}$ have at most two real zeros (see [10], page 23), it follows that $H(0, x, 1, \beta) > 0$ for $0 < x < 1$. For the same reason, the function $x \mapsto H(x, x_2, x_3, \beta)$ vanishes only at $x = x_2$ and $x = x_3$. Hence, $\text{sign}\, H(x_1, x_2, x_3, \beta) = \text{sign}\, H(0, x_2, x_3, \beta)$ if the ordering (3.11) holds.

Similarly, it follows that $\operatorname{sign} H(0, x_2, x_3, \beta) = \operatorname{sign} H(0, x_2, 1, \beta) = +1$. Recalling the statement on the coefficients $a_1$, $a_2$ and $a_3$ in (3.10), we see that the summands with $k = 1$ and $k = 3$ diminish the sum and it follows that

$$I_3(x_1, x_2, x_3, \beta) \leq I_3(0, x_2, 1, \beta).$$

This not only yields the cited result regarding the location of the smallest and largest support point of local $D$-optimal designs, but additionally provides the bound

$$I_3(x_1, x_2, x_3, \beta) \leq \sum_{k=1}^{3} I_3(0, x_k, 1, \beta),$$

which holds for support points $x_1$, $x_2$, $x_3$ in any order. In particular, we have a bound with a decomposition as stated in (2.9). Therefore, the exponential regression model can now be treated by arguments similar to those given in the previous examples, although the technical details are more involved. In particular, we have from (3.9), for large $\beta$,

$$I_3(0, x, 1, \beta) \leq \frac{2}{3} x^2 e^{-2\beta x}, \qquad |M(\xi[\beta], \beta)| = \frac{1}{3^3} \sup_x I_3(0, x, 1, \beta) \geq \frac{1}{27} \frac{1}{3e^2 \beta^2}.$$

Although there is no simple representation for the determinant $|M(\xi[\beta], \beta)|$, the estimates above are comparable to the corresponding equations (3.3) and (3.4) for the one-dimensional exponential model considered in Example 3.2 (the estimates differ only by constants). Thus, we conclude, similarly as in the previous examples, that the number of support points of the standardized maximin $D$-optimal design in the model (3.8) is unbounded for sufficiently large $\beta_{\max}$.

REMARK 3.5. At first glance, the results of Theorem 3.1 and the examples (including those in the technical report of Braess and Dette [1]) are surprising because it has never been observed in numerical studies that the number of support points of the standardized maximin $D$-optimal design substantially exceeds the number of parameters. To our knowledge, the numerical results of Example 3.2 are the first ones in this direction. However, it follows from the proof of Theorem 3.1 in the Appendix that the construction of a design with more than $N$ support points outperforming a given design requires a very large parameter space in the maximin $D$-optimality criterion. In particular, the number of support points of the standardized maximin $D$-optimal design may increase so slowly that it is almost impossible to decide by means of numerical computation whether it is asymptotically bounded or not. Thus, in practice, optimal designs with a large number of support points will only be observed if a very large parameter space is involved.

REMARK 3.6. It was pointed out by a referee that the minimum in the standardized maximin criterion is usually calculated over a finite grid in the interval $[\beta_{\min}, \beta_{\max}]$. A careful inspection of the proof of Theorem 3.1 shows that the results of this section can be extended to such situations.

**4. Bayesian $D$-optimal designs.** We now turn to analogous questions for the Bayesian $D$-optimality criterion (2.1). When Bayesian $D$-optimal designs are considered, it does not make a difference whether the information matrix or its standardized analogue is considered. The difference between the criterion

$$
\begin{aligned}
\Psi_{\rm st}(\xi) &= \int_{\mathcal{B}} \log \frac{|M(\xi,\beta)|}{|M(\xi[\beta],\beta)|} \pi(d\beta) \\
&= \int_{\mathcal{B}} [\log|M(\xi,\beta)| - \log|M(\xi[\beta],\beta)|]\pi(d\beta)
\end{aligned}
\tag{4.1}
$$

and the function in (2.1) is a constant that does not depend on the design $\xi$. In this case, uncertainty can be directly specified by a prior, which is supported on the interval $\mathcal{B} = [\beta_{\min}, \beta_{\max}]$, where $-\infty \leq \beta_{\min} \leq \beta_{\max} \leq \infty$. For example, one might increase the support of the prior without changing its shape or one could fix the support $\mathcal{B}$ of $\pi$ and change the shape such that its variance increases. The following result covers both cases.

THEOREM 4.1. *Assume that (2.9) holds and that $Q$ is uniformly decreasing with respect to the scale $\ell$ in the sense of Definition 2.1, where the function $\varphi$ satisfies*

$$\varphi(z) \leq c_1 e^{-|z|^\gamma} \tag{4.2}$$

*for some positive constants $c_1$, $\gamma$, and that the prior and the function $\ell$ in Definition 2.1 satisfy, for all measurable sets $B \subset \mathcal{B} = [\beta_{\min}, \beta_{\max}]$,*

$$\int_B \frac{c_3}{\ell(\mathcal{B})} \ell(d\beta) \leq \int_B \pi(d\beta) \tag{4.3}$$

*for some positive constant $c_3$. Let $N \in \mathbb{N}$ be given. If $\ell(\beta_{\max}) - \ell(\beta_{\min})$ is sufficiently large, then the Bayesian $D$-optimal design with respect to the prior $\pi$ on the interval $\mathcal{B}$ is supported at more than $N$ points.*

REMARK 4.2. Note that increasing the interval $\mathcal{B}$ in the optimality criterion (2.1) such that condition (4.3) is satisfied also changes the prior $\pi$ on $\mathcal{B}$. A typical example is the uniform distribution on the set $\mathcal{B}$ for which the assumption (4.3) is obviously satisfied if $\ell(\beta) = \beta$ or $\ell(\beta) = \log \beta$. In this case, the shape of the prior does not change, as will be illustrated in Example 4.3. On the other hand, uncertainty can also be quantified by changing the shape of the prior and, in this case, the function $\ell$ usually changes with $\pi$ (see Examples 4.4 and 4.5 below).

TABLE 2
*Bayesian D-optimal designs with respect to a uniform distribution on the interval $[1, B]$ in the exponential regression model (3.2). First row: support points; second row: weights*

| B     | 10    | 40    | 50    | 100   | 200   | 300   | 3000   |
|-------|-------|-------|-------|-------|-------|-------|--------|
| $x_1$ | 0.182 | 0.048 | 0.038 | 0.019 | 0.010 | 0.006 | 0.0006 |
| $w_1$ | 1.000 | 0.981 | 0.973 | 0.962 | 0.959 | 0.957 | 0.951  |
| $x_1$ |       | 0.354 | 0.318 | 0.215 | 0.134 | 0.084 | 0.009  |
| $w_2$ |       | 0.019 | 0.027 | 0.038 | 0.041 | 0.037 | 0.039  |
| $x_3$ |       |       |       |       |       | 0.236 | 0.055  |
| $w_3$ |       |       |       |       |       | 0.006 | 0.006  |
| $x_4$ |       |       |       |       |       |       | 1.000  |
| $w_4$ |       |       |       |       |       |       | 0.004  |

EXAMPLE 4.3. Consider the exponential regression model (3.2) of Example 3.2. Obviously, the function $\varphi$ in (3.6) also satisfies the stronger assumptions in Theorem 4.1. As a consequence, the number of support points of Bayesian $D$-optimal designs with respect to a uniform distribution is unbounded if the support of the prior is increased. Table 2 shows the Bayesian $D$-optimal designs corresponding to the situation considered in Table 1. Note that the standardized maximin $D$-optimal designs have remarkably more support points than the Bayesian $D$-optimal designs with respect to the uniform prior.

EXAMPLE 4.4. Consider the logistic regression model $Y \sim \text{Bin}(1, \eta(x, \beta))$ with

(4.4) $$\eta(x,\beta) = \frac{1}{1 + e^{x-\beta}}, \qquad x \in [0, \infty), \beta \in [0, \infty),$$

where the Fisher information of the parameter $\beta$ at point $x$ is given by

$$I(x, \beta) = \frac{e^{x-\beta}}{(1 + e^{x-\beta})^2}.$$

For any $a \in (0, 1)$, we consider the prior

$$\pi_a(d\beta) = cae^{-a\beta} I_{[0, 1/a)}(\beta) \, d\beta,$$

where $c = (1 - e^{-1})^{-1}$. Note that the expectation value and the variance of $\pi_a$ are proportional to $1/a$ and $1/a^2$, respectively. If we define

$$\ell_a(d\beta) = ca^{1/2} e^{-a\beta} I_{[0, 1/a)}(\beta) \, d\beta,$$

then $\ell_a(\mathcal{B}) = a^{-1/2} \to \infty$ if $a \to 0$, and a straightforward calculation shows that condition (4.3) holds with $c_3 = 1$.

The local $D$-optimal design is a one-point design concentrating its mass at the point $x[\beta] = \beta$ with $M(\xi[\beta], \beta) = 1/4$. Hence,

$$Q(\beta, \tilde{\beta}) = \frac{4e^{\tilde{\beta}-\beta}}{(1+e^{\tilde{\beta}-\beta})^2} \leq 4e^{|\beta-\tilde{\beta}|}.$$

An application of the mean value theorem shows that condition (2.4) is satisfied with $\varphi(z) = 4e^{-|z|}$ if $a \leq c^{-2}$. Moreover, $Q(\beta, \tilde{\beta}) \geq \frac{1}{2}$ if $|\ell_a(\beta) - \ell_a(\tilde{\beta})| \leq a^{1/2}$ and Theorem 4.1 applies. If $a \to 0$, the quantity $\ell_a(\mathcal{B}) = a^{-1/2}$ is sufficiently large and the number of support points of the Bayesian $D$-optimal design with respect to the prior $\pi_a$ in the logistic regression model (4.4) exceeds any given bound $N \in \mathbb{N}$.

EXAMPLE 4.5. As pointed out by a referee, it is worthwhile to mention that Theorem 4.1 also applies to discrete priors (where its proof has to be slightly modified). Consider, for example, the logistic regression model (4.4) and a uniform prior $\pi_L$ on the set $\mathcal{M}_L = \{1, \ldots, L\}$. If $\ell$ is the distribution function of the discrete measure with mass 1 at each element of $\mathcal{M}_L$, we have $\ell(\mathcal{B}) = L$ and there is equality in (4.3) with $c_3 = 1$ (note that $\ell$ is a step function with jumps of size 1 at each element of $\mathcal{M}_L$). Consequently, we have, for all $\beta, \tilde{\beta} \in \mathrm{supp}(\pi)$,

$$Q(\beta, \tilde{\beta}) = \frac{4e^{\tilde{\beta}-\beta}}{(1+e^{\tilde{\beta}-\beta})^2} \geq \frac{4e^{-1}}{(1+e^{-1})^2} \geq \frac{1}{2}$$

if $|\ell(\beta) - \ell(\tilde{\beta})| \leq 1$. Moreover, $Q(\beta, \tilde{\beta}) \leq 4ee^{-|\lfloor \beta \rfloor - \lfloor \tilde{\beta} \rfloor|}$ and it follows from Theorem 4.1 that the number of support points of the Bayesian $D$-optimal design with respect to the prior $\pi_L$ becomes arbitrarily large as $L \to \infty$. We finally mention a consequence of Carathéodory's theorem. For a discrete prior with $L$ support points, there exists a Bayesian $D$-optimal design with at most $Lm(m+1)/2$ support points; see [12]. This bound reduces to $L$ and converges to $\infty$ in the present case.

REMARK 4.6. Note that Theorem 4.1 requires stronger decay of the function $\varphi$ than Theorem 3.1. As a consequence, the number of support points of Bayesian $D$-optimal designs usually increases more slowly with the length of the parameter space compared to the maximin case. We have illustrated this fact in Examples 3.2 and 4.3, where we compared the standardized maximin and the Bayesian $D$-optimal design with respect to the uniform distribution. On the other hand, it follows from the proof of Theorem 4.1 in the Appendix that a similar result holds for the Bayesian $A$-optimality criterion

$$\int_{\mathcal{B}} \frac{\operatorname{trace}(M^{-1}(\xi, \beta))}{\operatorname{trace}(M^{-1}(\xi[\beta], \beta))} \pi(d\beta),$$

where condition (4.2) can be replaced by the weaker condition (3.1). This explains the empirical results of [2] that Bayesian $D$-optimal designs usually have more support points than Bayesian $A$-optimal designs.

**5. Conclusions.** When efficient designs in nonlinear regression models are constructed, it has been observed numerically by many authors that the number of support points of Bayesian and maximin $D$-optimal designs increases with the amount of uncertainty about a priori knowledge of the location of the nonlinear parameters. In this paper, we have established sufficient conditions under which the number of support points of Bayesian and maximin $D$-optimal designs can become arbitrarily large if the prior information on the unknown nonlinear parameters is diminished. The essential condition is the decay of the efficiency for large deviations between the specified and "true" parameter. The conditions apply to many of the commonly used regression models. In fact, we did not find any model where these conditions are not satisfied.

For the sake of brevity and a clear presentation, we have restricted our investigations to nonlinear models where one parameter appears nonlinearly in the Fisher information. However, our approach can also be applied to models with more nonlinear parameters, although some of the arguments have to be adapted. The main idea is to introduce an appropriate norm for high-dimensional nonlinear parameters which generalizes the distance (2.6). These arguments are outlined in Appendix A.3. A similar result is also available for the Bayesian $D$-optimality criterion by combining this argument with the results of Section 4. Moreover, Theorems 3.1 and 4.1 can also be extended to nonrectangular regions.

In this paper, we have made a general statement on the structure of optimal designs with respect to the standardized maximin and Bayesian $D$-optimality criteria, which is important for a better understanding of these sophisticated optimality criteria. For a given model of interest, our methodology can be used to prove a phenomenon which was conjectured for a long time in the literature. In all examples that we have investigated, the developed theory was applicable and we were able to prove that the number of support points of the standardized maximin and Bayesian $D$-optimal designs exceeds any given bound if the knowledge about the underlying parameter space is diminished. Moreover, we have also provided some explanation as to why standardized maximin $D$- and Bayesian $A$-optimal designs usually have more support points than Bayesian $D$-optimal designs.

We finally mention that the results for the Bayesian $D$-optimality criterion will have applications for estimating mixture distributions. To be precise, it was pointed out by Lindsay [14] that the determination of the ML-estimate of a mixture distribution corresponds to a Bayesian $D$-optimal design problem in a one-parameter nonlinear model. It therefore follows from the results of the present paper that in many models, the number of components of the estimated mixture distribution increases with the sample size.

## APPENDIX: PROOFS

**A.1. Proof of Theorem 3.1.** The proof consists of two steps. Set $B = \ell(\beta_{\max}) - \ell(\beta_{\min})$. First we show that, for an arbitrary design, say $\xi_N$, with $N$ support points, it follows that

$$(A.1) \quad \Phi(\xi_N) = \min\left\{ \frac{|M(\xi_N, \beta)|}{|M(\xi[\beta], \beta)|} \Big| \beta \in [\beta_{\min}, \beta_{\max}] \right\} \leq d_1(N+1)B^{-\gamma},$$

where $d_1$ is a positive constant not depending on $B$ and $\gamma > m$. Second we show that there exists a design $\xi_n$ (with at least $n$ support points) on $\mathcal{X}$ such that

$$(A.2) \quad \Phi(\xi_n) \geq \frac{d_2}{B^m}$$

for some positive constant $d_2$ not depending on $B$. Since $\gamma > m$, given $N$, we have

$$d_1(N+1)B^{-\gamma} < \frac{d_2}{B^m}$$

if $B$ is sufficiently large, and the optimal design is supported at more than $N$ points in this case. This proves the assertion. For the sake of a transparent representation, we begin with a proof of the estimates (A.1) and (A.2) in the case $m = 1$. The general case will be treated in a second step (B), while we prove in part (C) the remaining assertion of Theorem 3.1, considering the case where the local $D$-optimal designs are minimally supported.

(A) *The case $m = 1$*: To verify the estimate (A.1), let $\xi_N = \sum_{k=1}^{N} w_k \delta_{x_k}$ denote any design with mass $w_k$ at the point $x_k$ ($k = 1, \ldots, N$). Here, $\delta_{x_k}$ denotes the Dirac measure at the point $x_k$. Then

$$M(\xi_N, \beta) = \sum_{k=1}^{N} w_k I(x_k, \beta).$$

By assumption (2.9), there exist real numbers $\beta_{\min} \leq \beta_1 < \cdots < \beta_N \leq \beta_{\max}$ such that the inequality

$$(A.3) \quad M(\xi_N, \beta) \leq \sum_{k=1}^{N} w_k M(\xi[\beta_k], \beta) = M(\xi[\beta], \beta) \sum_{k=1}^{N} w_k Q(\beta, \beta_k)$$

holds for all $\beta \in \mathcal{B}$. For convenience, we put $\beta_0 = \beta_{\min}$, $\beta_{N+1} = \beta_{\max}$. Now, at least one gap between the numbers $\ell(\beta_k)$ must be large. Specifically, there exists an index $j \in \{0, \ldots, N\}$ such that

$$\text{(A.4)} \qquad \ell(\beta_{j+1}) - \ell(\beta_j) \geq \frac{\ell(\beta_{N+1}) - \ell(\beta_0)}{N+1} = \frac{B}{N+1}.$$

We consider the inequality at the point $\bar{\beta}$ defined by $\ell(\bar{\beta}) = \frac{1}{2}[\ell(\beta_j) + \ell(\beta_{j+1})]$ and derive from (A.4),

$$\text{(A.5)} \qquad |\ell(\bar{\beta}) - \ell(\beta_k)| \geq \frac{1}{2}(\ell(\beta_{j+1}) - \ell(\beta_j)) \geq \frac{B}{2(N+1)}$$

for all $k \in \{0, 1, 2, \ldots, N+1\}$. We now use inequality (A.3) and the definition of $Q$ in (2.3), and obtain from assumption (2.4), (3.1) and (A.5),

$$\text{(A.6)} \qquad \begin{aligned} M(\xi_N, \bar{\beta}) &\leq \sum_{k=1}^{N} w_k Q(\bar{\beta}, \beta_k) M(\xi[\bar{\beta}], \bar{\beta}) \\ &\leq \sum_{k=1}^{N} w_k \varphi(\ell(\bar{\beta}) - \ell(\beta_k)) M(\xi[\bar{\beta}], \bar{\beta}) \\ &\leq c_1 \left(\frac{B}{2(N+1)}\right)^{-\gamma} M(\xi[\bar{\beta}], \bar{\beta}) \\ &= \frac{c_1 (2N+2)^\gamma}{B^\gamma} M(\xi[\bar{\beta}], \bar{\beta}) \end{aligned}$$

for some positive constant $c_1$. We set $d_1 = c_1(2N+2)^\gamma$ and the proof of (A.1) is complete.

For proof of the lower bound (A.2), we may restrict ourselves to the case $B \geq 4\lambda$, where $\lambda$ is the constant defined in Definition 2.1(ii). We choose $n = \lceil \frac{1}{2} B/\lambda \rceil \leq B/\lambda$ and fix $\beta_k$ such that

$$\text{(A.7)} \qquad \ell(\beta_k) = \ell(\beta_{\min}) + (2k-1)\frac{B}{2n} \qquad (k = 1, 2, \ldots, n).$$

Note that these points are contained in the interval $[\beta_{\min}, \beta_{\max}]$ and that their distance is at most $2\lambda$. Let $\xi[\beta_k]$, $k = 1, \ldots, n$, denote the corresponding local $D$-optimal designs and define

$$\text{(A.8)} \qquad \xi_n = \frac{1}{n} \sum_{k=1}^{n} \xi[\beta_k].$$

The design $\xi_n$ has at least $n$ support points and its information matrix satisfies

$$\text{(A.9)} \qquad M(\xi_n, \beta) = \frac{1}{n} \sum_{k=1}^{n} M(\xi[\beta_k], \beta).$$

Obviously, given $\beta \in [\beta_{\min}, \beta_{\max}]$, there exists an index $j = j_\beta$ such that
$$|\ell(\beta) - \ell(\beta_j)| \leq \lambda.$$
By construction, $M(\xi[\beta_j], \beta) = Q(\beta, \beta_j) M(\xi[\beta], \beta) \geq \frac{1}{2} M(\xi[\beta], \beta)$. Since all terms in the sum (A.9) are nonnegative, it follows that for all $\beta \in \mathcal{B}$,
$$M(\xi_n, \beta) \geq \frac{1}{n} M(\xi[\beta_j], \beta) \geq \frac{1}{2n} M(\xi[\beta], \beta) \geq \frac{\lambda}{B} M(\xi[\beta], \beta).$$
Recalling the definition of the standardized maximin criterion in (2.2), we conclude that $\Phi(\xi_n) \geq \lambda/2B$. With the choice $d_2 = \lambda/2$, we have proven the lower bound (A.2), and the proof in the case $m = 1$ is complete.

(B) *The case $m \geq 1$:* Let $\xi_N$ be a design with masses $w_k$ at the points $x_k \in \mathcal{X}$ ($k = 1, \ldots, N$), and for any tuple $(i_1, \ldots, i_m)$ with $1 \leq i_1 < \cdots < i_m \leq N$, let $\xi[\beta^{(1)}_{i_1,\ldots,i_m}], \ldots, \xi[\beta^{(m)}_{i_1,\ldots,i_m}]$ denote the designs corresponding to the points $x_{i_1}, \ldots, x_{i_m}$ by inequality (2.9). Using the definition of $Q$ in (2.3) and the Cauchy–Binet formula, we obtain

$$|M(\xi_N, \beta)| = \sum_{1 \leq i_1 < \cdots < i_m \leq N} w_{i_1} \cdots w_{i_m} I_m(x_{i_1}, \ldots, x_{i_m}, \beta)$$

$$\text{(A.10)} \qquad \leq c_0 \sum_{1 \leq i_1 < \cdots < i_m \leq N} w_{i_1} \cdots w_{i_m} \sum_{j=1}^m |M(\xi[\beta^{(j)}_{i_1,\ldots,i_m}], \beta)|$$

$$= c_0 |M(\xi[\beta], \beta)| \sum_{j=1}^m \sum_{1 \leq i_1 < \cdots < i_m \leq N} w_{i_1} \cdots w_{i_m} Q(\beta, \beta^{(j)}_{i_1,\ldots,i_m}).$$

Note that there are $m\binom{N}{m}$ terms in this sum and that inequality (A.10) corresponds to (A.3) in the proof of the case $m = 1$. Therefore, ordering the points $\beta^{(j)}_{i_1,\ldots,i_m}$ and using exactly the same arguments as in the proof of part (A) yields the upper bound $\Phi(\xi_N) \leq d_1 B^{-\gamma}$ for some constant $d_1$ and $\gamma > m$.

In order to prove the corresponding lower bound, we define $n = \lceil \frac{1}{2} B/\lambda \rceil$ and again consider the quantities $\beta_k$ defined by (A.7). Let $\xi[\beta_k]$ denote the corresponding local $D$-optimal design, define the design $\xi_n$ by (A.8) and denote by $\tilde{x}_i$ and $\tilde{w}_i$ the corresponding support points and weights of $\xi_n$, respectively. For any $\beta$, there exists a $\beta_j$ such that $|\ell(\beta) - \ell(\beta_j)| \leq \lambda$ and we denote the support points and weights of the corresponding local $D$-optimal design $\xi[\beta_j]$ by $x_i$ and $w_i$, respectively. By the Cauchy–Binet formula, we obtain

$$|M(\xi_n, \beta)| = \sum_{i_1 < \cdots < i_m} \tilde{w}_{i_1} \cdots \tilde{w}_{i_m} I_m(\tilde{x}_{i_1}, \ldots, \tilde{x}_{i_m}, \beta)$$

$$\text{(A.11)} \qquad \geq \frac{1}{n^m} \sum_{i_1 < \cdots < i_m} w_{i_1} \cdots w_{i_m} I_m(x_{i_1}, \ldots, x_{i_m}, \beta)$$

$$= \frac{1}{n^m} |M(\xi[\beta_j], \beta)|,$$

where the last inequality follows by omitting all terms containing points which are not in the support of the local $D$-optimal design $\xi[\beta_j]$. Using assumption (2.5), we therefore obtain

$$\text{(A.12)} \quad |M(\xi_n, \beta)| \geq \frac{1}{n^m}|M(\xi[\beta], \beta)|Q(\beta, \beta_j) \geq \frac{1}{2n^m}|M(\xi[\beta], \beta)|$$

and the same argument as presented in the proof for the case $m = 1$ shows the lower bound $\Phi(\xi_n) \geq cd_2/B^m$ for some positive constant $d_2$. The assertion in the case $m \geq 1$ now follows by the same arguments as given in the first part of the proof.

(C) *Proof of the remaining assertion.* If the local $D$-optimal designs are supported at $m$ points, the corresponding weights all equal $1/m$ (see [18]) and the estimate in (A.11) can be improved as follows. Let $\bar{x}_1, \ldots, \bar{x}_{m_\eta}$ denote the points of the design $\xi_n$, which are support points of any local $D$-optimal design, and define $\bar{w}_1, \ldots, \bar{w}_{m_\eta}$ to be the corresponding weights. Note that $\bar{w}_i = 1/m$, $i = 1, \ldots, m_\eta$; then

$$|M(\xi_n, \beta)|$$
$$\geq \sum_{i_{m_\eta+1} < \cdots < i_m} \bar{w}_1 \cdots \bar{w}_{i_{m_\eta}} \tilde{w}_{i_{m_\eta+1}} \cdots \tilde{w}_{i_m} I_m(\bar{x}_1, \ldots, \bar{x}_{m_\eta}, \tilde{x}_{i_{m_\eta+1}}, \ldots, \tilde{x}_{i_m}, \beta)$$
$$\geq \frac{1}{m^m n^{m-m_\eta}} I_m(x_1, \ldots, x_m, \beta) = \frac{1}{n^{m-m_\eta}}|M(\xi[\beta_j], \beta)|.$$

Here, the first inequality follows by considering only the terms for which $I_m$ contains the common support points $\bar{x}_1, \ldots, \bar{x}_{m_\eta}$, while the second inequality is obtained by considering only the term corresponding to the local $D$-optimal design $\xi[\beta_j]$. The same argument as used for (A.12) now shows that $\Phi(\xi_n) \geq d_2/B^{m-m_\eta}$ for some positive constant $d_2$, and the assertion now follows, as explained in the first part of the proof.

**A.2. Proof of Theorem 4.1.** We restrict ourselves to the case $m = 1$. The case $m \geq 1$ can be obtained by adapting arguments from part (B) in the proof of Theorem 3.1.

Moreover, for convenience, in a first step, we assume that the given transformation $\ell$ is the identity and define $B = \ell(\mathcal{B}) = \beta_{\max} - \beta_{\min}$. For a given design $\xi_N$ with $N$ support points, we know that inequality (A.3) holds. With assumption (4.3) and the notation $\beta_0 = \beta_{\min}$, $\beta_{N+1} = \beta_{\max}$ and $\Delta_j = \frac{1}{2}(\beta_{j+1} - \beta_j)$, we estimate the contribution of the interval $[\beta_j, \beta_j + \Delta_j]$ to the Bayesian $D$-optimality criterion via

$$\int_{\beta_j}^{\beta_j + \Delta_j} \log \frac{|M(\xi_N, \beta)|}{|M(\xi[\beta], \beta)|} \pi(d\beta)$$

$$\leq \frac{c_3}{B} \int_{\beta_j}^{\beta_j+\Delta_j} \log \sum_{k=1}^{N} w_k \frac{|M(\xi[\beta_k],t)|}{|M(\xi[t],t)|} \, dt$$

$$\leq \frac{c_3}{B} \int_{\beta_j}^{\beta_j+\Delta_j} \log \sum_{k=1}^{N} w_k \varphi(t-\beta_k) \, dt$$

$$\leq \frac{c_3}{B} \int_{\beta_j}^{\beta_j+\Delta_j} \log(c_1 \exp(-|t-\beta_j|^\gamma)) \, dt$$

$$\leq \frac{c_3}{B} \int_0^{\Delta_j} [\log c_1 - z^\gamma] \, dz = \frac{c_3}{B} \left[ \Delta_j \log c_1 - \frac{1}{1+\gamma} \Delta_j^{1+\gamma} \right].$$

The same bound is derived for the interval $[\beta_j + \Delta_j, \beta_{j+1}]$. Summing over all intervals of this form, we conclude that

$$\Psi_{\text{st}}(\xi_N) \leq \frac{2c_3}{B} \sum_{j=0}^{N} \left[ \Delta_j \log c_1 - \frac{1}{1+\gamma} \Delta_j^{1+\gamma} \right].$$

Since $\sum_{j=0}^{N} \Delta_j = \frac{B}{2}$ and the function $z^{1+\gamma}$ is strictly convex, the right-hand side attains its maximum if all $\Delta_j$'s are equal, and we obtain the upper bound

(A.13) $$\Psi_{\text{st}}(\xi_N) \leq c_4 \left[ \log c_1 - \frac{B^\gamma}{(1+\gamma)} \right]$$

for some positive constant $c_4$. Note that the right-hand side of this inequality is dominated by the term with $B^\gamma$ when $B \to \infty$.

The construction of a better design with respect to the Bayesian optimality criterion follows the arguments given in part (A) of the proof of Theorem 3.1. Let $\lambda > 0$ be defined by

$$Q(\beta, \tilde{\beta}) \geq \tfrac{1}{2} \qquad \text{whenever } |\beta - \tilde{\beta}| \leq \lambda.$$

Set $n = \lceil \frac{B}{2\lambda} \rceil$ and $\beta_j = \beta_{\min} + (2j-1)\lambda$ for $k = 1, 2, \ldots, n$. We choose a design $\xi$ with (at least) $n$ support points such that

$$M(\xi, \beta) = \frac{1}{n} \sum_{k=1}^{n} M(\xi[\beta_k], \beta)$$

[see the identity in (A.9)]. For any given $\beta \in [\beta_{\min}, \beta_{\max}]$, there exists a $\beta_j$ with $|\beta - \beta_j| \leq \lambda$ satisfying (2.5). Therefore, it follows for all $\beta \in \mathcal{B}$ that

$$M(\xi, \beta) \geq \frac{1}{n} M(\xi[\beta_j], \beta) = \frac{1}{n} Q(\beta, \beta_j) M(\xi[\beta], \beta) \geq \frac{1}{2n} M(\xi[\beta], \beta).$$

Hence,

$$\Psi_{\text{st}}(\xi_n) \geq \int_{\mathcal{B}} \log \frac{1}{2n} \pi(d\beta) = \log \frac{1}{2n} \geq -\log B + \log \lambda.$$

This value is larger than the upper bound (A.13) if $B$ is sufficiently large. Therefore, a design with $N$ support points cannot be optimal if $B$ is sufficiently large. Thus far, we have restricted ourselves to the case $\ell(\beta) = \beta$. The general case proceeds in exactly the same way, where one must choose $dt = \ell(d\beta)/\ell(\mathcal{B})$ in the first integral of the proof, and the boundaries of the intervals must be adapted. The details are left to the reader.

**A.3. A comment on more "nonlinear" parameters.** In this paragraph, we briefly describe how the arguments need to be changed if there exist $p > 1$ nonlinear parameters, say $\beta = (\beta_1, \ldots, \beta_p)$, which appear nonlinearly in the Fisher information matrix. For the sake of brevity, we consider the standardized maximin criterion with a $p$-dimensional cube $\mathcal{B}$. First, note that condition (2.9) does not depend on the dimension of the parameter $\beta$. Second, let $d$ denote a norm on $\mathbb{R}^p$ and replace conditions (2.4) and (2.5) by

$$(A.14) \qquad Q(\beta, \tilde{\beta}) \leq \varphi(d(\beta, \tilde{\beta}))$$

and

$$(A.15) \qquad Q(\beta, \tilde{\beta}) \geq \tfrac{1}{2} \qquad \text{whenever } d(\beta, \tilde{\beta}) \leq \lambda,$$

respectively. We then show that the number of support points of the standardized maximin $D$-optimal design becomes arbitrarily large if the volume of the cube $\mathcal{B}$ converges to infinity. In other words, Theorem 3.1 also holds in the case where the $p > 1$ parameters appear nonlinearly in the Fisher information.

For this, we note that the proof of Theorem 3.1 is performed by establishing the bounds (A.1) and (A.2). Let $\xi_N$ denote the design considered in (A.10) and note that the estimate (A.10) does not depend on the dimension of $\beta$. For $r > 0$, define the ball with center $\beta$ and radius $r$ by

$$U_r(j, i_1, \ldots, i_m) = \{x \in \mathbb{R}^p | d(x, \beta^{(j)}_{i_1, \ldots, i_m}) \leq r\}.$$

There exists a minimal $r_{\min}$ such that $\mathcal{B}$ can be covered by balls of this type, that is,

$$\mathcal{B} \subset \bigcup_{j=1}^m \bigcup_{i_1 < \cdots < i_m} U_{r_{\min}}(j, i_1, \ldots, i_m).$$

Obviously, we have, for some constant $c > 0$, that $r_{\min} > cB$, where $B$ denotes the $p$th root of the volume of $\mathcal{B}$. Consequently, there exists a $\tilde{\beta} \in \mathcal{B}$ such that

$$d(\tilde{\beta}, \beta^{(j)}_{i;1,\ldots,i_m}) \geq r_{\min}/2 \geq cB.$$

Thus, replacing condition (2.4) by (A.14) in the argument (A.6) yields the upper bound (A.1) for any $N$-point design. The remaining inequality (A.2) is similarly obtained by covering $\mathcal{B}$ with balls of radius $\lambda$.

**Acknowledgments.** Parts of this paper were written during the second author's visit to the Institute of Statistics in Louvain-la-Neuve in September 2004 and this author would like to thank the Institute for its hospitality. The authors are also grateful to A. Pepelyshev for computational assistance and to Isolde Gottschlich who typed parts of this paper with considerable technical expertise. We are also grateful to two unknown referees for their constructive comments which led to a substantial improvement of the presentation.

## REFERENCES

[1] BRAESS, D. and DETTE, H. (2005). On the number of support points of maximin and Bayesian $D$-optimal designs in nonlinear regression models. Technical report. Available at http://www.ruhr-uni-bochum.de/mathematik3/preprint.htm.

[2] CHALONER, K. and LARNTZ, K. (1989). Optimal Bayesian design applied to logistic regression experiments. *J. Statist. Plann. Inference* **21** 191–208. MR0985457

[3] DETTE, H. (1997). Designing experiments with respect to "standardized" optimality criteria. *J. Roy. Statist. Soc. Ser. B* **59** 97–110. MR1436556

[4] DETTE, H. and BIEDERMANN, S. (2003). Robust and efficient designs for the Michaelis–Menten model. *J. Amer. Statist. Assoc.* **98** 679–686. MR2011681

[5] DETTE, H. and NEUGEBAUER, H.-M. (1996). Bayesian optimal one point designs for one parameter nonlinear models. *J. Statist. Plann. Inference* **52** 17–31. MR1391681

[6] DETTE, H. and NEUGEBAUER, H.-M. (1997). Bayesian $D$-optimal designs for exponential regression models. *J. Statist. Plann. Inference* **60** 331–349. MR1456635

[7] FORD, I., TORSNEY, B. and WU, C.-F. J. (1992). The use of a canonical form in the construction of locally optimal designs for non-linear problems. *J. Roy. Statist. Soc. Ser. B* **54** 569–583. MR1160483

[8] HAN, C. and CHALONER, K. (2003). $D$- and $c$-optimal designs for exponential regression models used in viral dynamics and other applications. *J. Statist. Plann. Inference* **115** 585–601. MR1985885

[9] IMHOF, L. A. (2001). Maximin designs for exponential growth models and heteroscedastic polynomial models. *Ann. Statist.* **29** 561–576. MR1863970

[10] KARLIN, S. and STUDDEN, W. J. (1966). *Tchebycheff Systems. With Applications in Analysis and Statistics*. Wiley, New York. MR0204922

[11] KIEFER, J. C. (1974). General equivalence theory for optimum designs (approximate theory). *Ann. Statist.* **2** 849–879. MR0356386

[12] LÄUTER, E. (1974). Methode der Planung eines Experiments für den Fall nichtlinearer Parametrisierung. *Math. Operationsforsch. Stat.* **5** 625–636. (In Russian.)

[13] LIEBIG, H.-P. (1988). Temperature integration by kohlrabi growth. *Acta Horticulturae* **230** 371–380.

[14] LINDSAY, B. G. (1983). The geometry of mixture likelihoods: A general theory. *Ann. Statist.* **11** 86–94. MR0684866

[15] MUKHOPADHYAY, S. and HAINES, L. M. (1995). Bayesian $D$-optimal designs for the exponential growth model. *J. Statist. Plann. Inference* **44** 385–397. MR1332682

[16] MÜLLER, CH. H. (1995). Maximin efficient designs for estimating nonlinear aspects in linear models. *J. Statist. Plann. Inference* **44** 117–132. MR1323074

[17] PRONZATO, L. and WALTER, E. (1985). Robust experiment design via stochastic approximation. *Math. Biosci.* **75** 103–120. MR0800967

[18] PUKELSHEIM, F. (1993). *Optimal Design of Experiments*. Wiley, New York. MR1211416
[19] WONG, W. K. (1992). A unified approach to the construction of minimax designs. *Biometrika* **79** 611–619. MR1187611

RUHR-UNIVERSITÄT BOCHUM
FAKULTÄT FÜR MATHEMATIK
44780 BOCHUM
GERMANY
E-MAIL: dietrich.braess@rub.de
       holger.dette@rub.de